\documentclass[12pt,english]{article}
\usepackage{fullpage}
\usepackage[latin1]{inputenc}
\usepackage{babel}
\usepackage{exscale,xspace}
\usepackage{amsfonts,amssymb}

\usepackage[
  pdfauthor={Bernd C. Kellner},
  pdfkeywords={Bernoulli numbers, Stirling numbers, Carmichael numbers,
    conjecture of Giuga, conjecture of Agoh},
  pdftitle={The Equivalence of Giuga's and Agoh's Conjectures},
  colorlinks,bookmarksnumbered]{hyperref}

\newcommand{\qed}{\parbox{0cm}{}\hspace*{\fill} $\Box$}

\newcommand{\ds}{\scriptstyle}
\newcommand{\binom}[2]{ {#1 \choose #2} }
\newcommand{\pdiv}{\mid}

\newcommand{\notdiv}{\nmid}

\newcommand{\mod}[1]{({\rm mod\ } #1)}
\newcommand{\valueat}[1]{{\ }_{\big| #1}}

\def\denom{\mathop{\rm denom}\nolimits}

\newcommand{\ZZ}{\mathbb{Z}}
\newcommand{\NN}{\mathbb{N}}

\newcommand{\refeqn}[1]{(\ref{#1})}

\newcommand{\stz}[2]{ \bigg\{ {#1 \atop #2} \bigg\} }
\newcommand{\ste}[2]{ \bigg[ {#1 \atop #2} \bigg] }
\newcommand{\ST}[2]{ \bigg< {#1 \atop #2} \bigg> }
\newcommand{\STL}[2]{ \big< {#1 \atop #2} \big> }
\newcommand{\stirling}{{\bf S}}

\renewcommand{\theequation}{\arabic{section}.\arabic{equation}}
\newcommand{\settheequation}[1]{\renewcommand{\theequation}{#1}}
\newcommand{\restoretheequation}{\renewcommand{\theequation}%
{\arabic{section}.\arabic{equation}}\addtocounter{equation}{-1}}

\newtheorem{prop}{Proposition}[section]
\newtheorem{theorem}[prop]{Theorem}
\newtheorem{corl}[prop]{Corollary}
\newtheorem{lemma}[prop]{Lemma}
\newtheorem{conj}[prop]{Conjecture}

\newtheorem{thdefin}[prop]{Definition}
\newtheorem{thremark}[prop]{Remark}
\newtheorem{thtable}[prop]{Table}

\newenvironment{defin}{\begin{thdefin}\rm}{\end{thdefin}}
\newenvironment{remark}{\begin{thremark}\rm}{\end{thremark}}

\newenvironment{proof}{\medskip\noindent%
\textsc{Proof.}}{\qed\medskip}
\newenvironment{proofof}[1]{\medskip\noindent%
\textsc{Proof of #1.}}{\qed\medskip}

\newenvironment{equationname}[1]{\settheequation{#1}\begin{equation}}%
{\end{equation}\restoretheequation}

\begin{document}

\title{The Equivalence of Giuga's and Agoh's Conjectures}
\author{Bernd C. Kellner}
\date{}
\maketitle

\abstract{In this paper we show the equivalence of the conjectures of Giuga and Agoh
in a direct way which leads to a combined conjecture. This conjecture is described
by a sum of fractions from which all conditions can be derived easily.}
\smallskip

\smallskip


\parindent 0cm

\section{Introduction}

The Bernoulli numbers $B_n$ are defined by the power series
\[
   \frac{z}{e^z-1} = \sum_{n=0}^\infty B_n \frac{z^n}{n!} \,,
     \qquad |z| < 2 \pi \,,
\]
where all numbers $B_n$ are zero with odd index $n > 1$.
First values are given by $B_0 = 1$, $B_1 = -\frac{1}{2}$,
$B_2 = \frac{1}{6}$, $B_4 = -\frac{1}{30}$.
The Bernoulli polynomials are similarly defined by the generating function
\[
   \frac{ze^{zx}}{e^z-1} = \sum_{n=0}^\infty B_n(x) \,\, \frac{z^n}{n!} \,,
      \qquad |z| < 2 \pi \,.
\]

One has the well-known properties $B_n(0) = B_n$, $B_n'(x) = n B_{n-1}(x)$, and
\[
   B_n(x) = \sum_{k=0}^n \binom{n}{k} B_k \, x^{n-k} \,.
\]

Define for positive integers $n$ and $m$ the summation formula of consecutive
integer powers by
\[
   S_n(m) = \sum_{k=1}^{m-1} k^n.
\]

\begin{theorem} \label{theor-sn-bn}
Let $n$ be a positive integer, then for real $x$
\[
   S_n(x) = \frac{1}{n+1} ( B_{n+1}(x) - B_{n+1} ) \,.
\]
\end{theorem}

These basic facts can be found in \cite[Chapter 15]{ireland90classical}.
Here we give a general congruence between $S_n$ and $B_n$
which is valid for arbitrary integers $m>1$ and even $n$.
We will prove it later using Stirling numbers
of the second kind. Throughout this paper, $p$ will denote a prime.

\begin{theorem} \label{theor-snm-mbn}
Let $n$, $m$ be positive integers with $m > 1$ and even $n$, then
\[
   S_n(m) \equiv m \, B_n \equiv
     - \!\! \sum_{\ds p \pdiv m \atop \ds p-1 \pdiv n} \!\! \frac{m}{p} \pmod{m} \,.
\]
\end{theorem}

The following basic lemma we will need later on.

\begin{lemma} \label{lem-divcongr}
Let $a$, $c$, and $m$ be positive integers with $a \pdiv m$, then
\[
   c \, \frac{m}{a} \equiv c' \, \frac{m}{a} \pmod{m} \qquad \mbox{for} \quad
     c \equiv c' \pmod{a} \,.
\]
\end{lemma}

\begin{proof}
Since $c \equiv c' \ \mod{a}$, there exists $k \in \ZZ$ with $c=ak+c'$. Hence
\[
   c \, \frac{m}{a} \equiv (ak+c') \, \frac{m}{a} \equiv
     km + c' \, \frac{m}{a} \equiv c' \, \frac{m}{a} \pmod{m} \,. \vspace*{-4ex}
\]
\end{proof}

\section{Equivalence of the conjectures}
\setcounter{equation}{0}

The following conjectures were independently formulated
by G.\,Giuga (ca 1950, see \cite{giuga50su}) and
by T.\,Agoh (ca 1990, see \cite{agoh95}).
In 1993, Agoh established the connection and equivalence of both conjectures.

\begin{conj}[Giuga, 1950]
Let $n$ be a positive integer with $n \geq 2$, then
\[
   S_{n-1}(n) \equiv -1 \pmod{n} \quad \Longleftrightarrow \quad n \ \mbox{is prime} \,.
\]
\end{conj}

\begin{conj}[Agoh, 1990] \label{conj-agoh}
Let $n$ be a positive integer with $n \geq 2$, then
\[
   n B_{n-1} \equiv -1 \pmod{n} \quad \Longleftrightarrow \quad n \ \mbox{is prime} \,.
\]
\end{conj}

The following theorem implies the equivalence of both conjectures.

\begin{theorem}
Let $n$ be a positive integer with $n \geq 2$, then
\[
   S_{n-1}(n) - n B_{n-1} \equiv \left\{
     \begin{array}{rl}
       n/2 \,, & (n \equiv 2 \ \mod{4}, \, n > 2) \\
       0   \,, & \mbox{otherwise}
     \end{array} \right.
     \pmod{n} \,.
\]
\end{theorem}

\begin{proof}
Case $n=2$ is trivial: $S_1(2) \equiv 2 B_1 \equiv 1 \ \mod{2}$.
For odd $n$ with $n \geq 3$, Theorem \ref{theor-snm-mbn} provides
\begin{eqnarray} \label{eqn-equiv-giuga-agoh}
   S_{n-1}(n) \equiv n \, B_{n-1} \equiv
     - \!\!\!\! \sum_{\ds p \pdiv n \atop \ds p-1 \pdiv n-1} \!\! \frac{n}{p} \pmod{n} \,.
\end{eqnarray}
For now, the cases $n \geq 4$ with even $n$ remain, then $B_{n-1} = 0$.
We only have to determine $S_{n-1}(n) \ \mod{n}$.
Since $n-1$ is odd, $\nu^{n-1} \equiv - (n-\nu)^{n-1} \ \mod{n}$ for $\nu = 1,\ldots,n/2$.
Hence, all elements of the sum cancel each other except $\nu=n/2$ in the middle.
Using Lemma \ref{lem-divcongr} provides
\[
   S_{n-1}(n) \equiv  \left( \frac{n}{2} \right)^{n-1}
   \equiv \frac{n}{2} \left( \frac{n}{2} \right)^{n-2}
   \equiv \left\{
     \begin{array}{rl}
       n/2 \,, &  (n \equiv 2 \ \mod{4}) \\
       0   \,, &  (n \equiv 0 \ \mod{4})
     \end{array} \right.
   \hspace*{-1.5ex} \pmod{n} \,,
\]
since $4 \pdiv n$ yields $n/2 \equiv 0 \ \mod{2}$ and
$n \equiv 2 \ \mod{4}$ yields $n/2 \equiv 1 \ \mod{2}$.
\end{proof}

Because both conjectures are equivalent in a simple manner,
it will be called, for now, the conjecture of Giuga-Agoh which can
be formulated in another way.
The congruence \refeqn{eqn-equiv-giuga-agoh} is valid
for odd $n$ and $n=2$. Moreover, we have for even $n$, $n \geq 4$
\[
   \sum_{\ds p \pdiv n \atop \ds p-1 \pdiv n-1} \!\! \frac{n}{p}
     \equiv \frac{n}{2} \not\equiv 1 \pmod{n} \,,
\]
since $p-1 \pdiv n-1$ is only valid for $p=2$.
In that case, we have $S_{n-1}(n) \not\equiv -1 \ \mod{n}$ and
$n B_{n-1} \equiv 0 \ \mod{n}$. Thereby, we obtain
another equivalent conjecture which
is described without Bernoulli numbers and summation function $S_n$.

\begin{conj}[Giuga-Agoh]
Let $n$ be a positive integer with $n \geq 2$, then
\begin{equationname}{G}  \label{eqn-giuga}
   \sum_{\ds p \pdiv n \atop \ds p-1 \pdiv n-1} \!\! \frac{n}{p} \equiv 1 \pmod{n}
     \quad \Longleftrightarrow \quad n \ \mbox{is prime} \,.
\end{equationname}
\end{conj}
\medskip

Each prime $p$ yields a trivial solution of \refeqn{eqn-giuga} with $n=p$.
Therefore any nontrivial solution $n$ of \refeqn{eqn-giuga} must be composite
and provides a counterexample of the conjecture of Giuga-Agoh.
So far, no counterexample was found.

\section{Conditions and properties}
\setcounter{equation}{0}

\begin{lemma} \label{lem-giuga-prop}
Let $n$ be a nontrivial solution of \refeqn{eqn-giuga}, then
\begin{enumerate}
\item $n$ is composite, odd, and squarefree.
\item $p \pdiv n/p-1$ for all prime divisors $p$ of $n$.
\item $p-1 \pdiv n-1$ and $p-1 \pdiv n/p-1$ for all prime divisors $p$ of $n$.
\end{enumerate}
\end{lemma}

\begin{proof} All properties will follow by congruence \refeqn{eqn-giuga}.
A nontrivial solution $n$ is composite. Let $p$, $q$ be prime divisors of $n$,
then one also has
\[
   \sum_{\ds q \pdiv n \atop \ds q-1 \pdiv n-1} \!\! \frac{n}{q} \equiv 1 \pmod{p} \,.
\]
In case $p \neq q$ the term $n/q \equiv 0 \ \mod{p}$ vanishes. Hence, only
$n/p \equiv 1 \ \mod{p}$ remains with condition $p-1 \pdiv n-1$, otherwise
the whole sum would vanish. (1), (2): $n/p \equiv 1 \ \mod{p}$ implies
$p^2 \notdiv n$ and $p \pdiv n/p-1$. Therefore $n$ is squarefree. Moreover
$n$ is odd, otherwise $p-1 \pdiv n-1$ is only valid for $n=2$.
(3): $p-1 \pdiv n-1$ yields $n \equiv 1 \ \mod{p-1}$ and $n/p \equiv 1/p \equiv 1 \ \mod{p-1}$.
\end{proof}

Giuga \cite{giuga50su} proved properties (1)-(3) and following statements.
Regarding properties (2) and (3) separately, one can achieve further
conditions and properties. Therefore we have to introduce some definitions
of Giuga and Carmichael \cite{carmichael10}.

\begin{defin}
A composite positive integer $m$ is called a \textsl{Carmichael number}, if
\[
   a^{m-1} \equiv 1 \pmod{m}
\]
is valid for all $(a,m)=1$.
\end{defin}

\begin{theorem}[Carmichael] \label{theor-carmichael-num}
A Carmichael number $m$ is odd, squarefree and has at least 3 prime factors.
Let $p$, $q$ be prime divisors of $m$, then the following conditions hold
\[
   p-1 \pdiv m-1 \,, \quad p-1 \pdiv m/p-1 \,, \quad q \not\equiv 1 \pmod{p} \,.
\]
\end{theorem}

The first three Carmichael numbers are $561 = 3 \cdot 11 \cdot 17$,
$1105 = 5 \cdot 13 \cdot 17$, and $1729 = 7 \cdot 13 \cdot 19$.
It was proven in \cite{alford93carmichael},
that infinitely many Carmichael numbers exist.

\begin{defin}
A composite positive integer $n$ is called a \textsl{Giuga number}, if
\begin{equation} \label{eqn-giuga-num}
   \sum_{p \pdiv n} \frac{1}{p} - \prod_{p \pdiv n} \frac{1}{p} \in \NN \,.
\end{equation}
\end{defin}
\smallskip

The previous definition was first given by Giuga.
The term \textsl{Giuga number} was introduced in \cite{borwein95survey} and \cite{borwein96giuga},
where one can also find the definition of a \textsl{Giuga sequence} as its generalization.
The first three Giuga numbers are $30 = 2 \cdot 3 \cdot 5$,
$858 = 2 \cdot 3 \cdot 11 \cdot 13$, and $1722 = 2 \cdot 3 \cdot 7 \cdot 41$.
Only even Giuga numbers have been found yet.
Now, properties of Lemma \ref{lem-giuga-prop} show that
a counterexample of \refeqn{eqn-giuga} must be both a Giuga number and
a Carmichael number.

\begin{theorem}
Let $n$ be a nontrivial solution of \refeqn{eqn-giuga}, then
\begin{enumerate}
\item $n$ is odd, squarefree, and has at least 9 prime factors.
\item $n$ is a Giuga number and $p \pdiv n/p-1$ for all prime divisors $p$ of $n$.
\item $n$ is a Carmichael number and $p-1 \pdiv n-1$ resp.\ $p-1 \pdiv n/p-1$ for all prime divisors $p$ of $n$.
\end{enumerate}
\end{theorem}

\pagebreak

\begin{proof}
We simply have to extend the results of Lemma \ref{lem-giuga-prop}.
Divide \refeqn{eqn-giuga} by $n$, then
\[
   \sum_{\ds p \pdiv n \atop \ds p-1 \pdiv n-1} \!\! \frac{1}{p} - \frac{1}{n}
     \equiv 0 \pmod{\ZZ} \,.
\]
(2): Since $n$ consists at least of two prime factors, the left side above must lie in $\NN$.
Without condition $p-1 \pdiv n-1$ this yields \refeqn{eqn-giuga-num} and $n$ must
be a Giuga number. (1): Let $p_\nu$ be the $\nu$-th prime. Then we have
$\sum_{\nu=2}^{9} 1/p_\nu < 1$. Since $n$ is odd, $n$ must have at least 9 prime factors.
(3): Theorem \ref{theor-carmichael-num} shows that conditions of (3)
identify $n$ as a Carmichael number.
\end{proof}

A Carmichael number has restrictions on its prime factors $p$, as seen
in Theorem \ref{theor-carmichael-num}. For $p \pdiv n$ there
cannot occur other prime factors $q$ of $n$ with $q \equiv 1 \ \mod{p}$.
On the other side, a Giuga number $n$ must satisfy
\[
    \sum_{p \pdiv n} \frac{1}{p} > 1 \,.
\]

In 1950, using these two properties Giuga
showed that a nontrivial solution $n$ respectively a counterexample
of \refeqn{eqn-giuga} must have more than 360 prime factors which
provides $n > 10^{1000}$. Bedocchi \cite{bedocchi85giuga} extended
this result to $n > 10^{1700}$ in 1985. Finally, in 1996,
D. Borwein, J. M. Borwein, P. B. Borwein, and Girgensohn \cite{borwein96giuga}
raised the limit to $n > 10^{13887}$ by further reducing of possible
cases. For corresponding methods see the references.
\smallskip

Agoh \cite{agoh95} also showed congruence \refeqn{eqn-giuga}
and \refeqn{eqn-equiv-giuga-agoh}, where the equivalence is
essentially derived by means of the Theorem of Clausen-von Staudt.
In \cite{agoh95} one additionally finds stronger conditions and
extended results on nontrivial solutions of \refeqn{eqn-giuga}.
Further examples of Giuga numbers are given in
\cite{borwein95survey} and \cite{borwein96giuga}, see also
\cite{butske99} for new solutions of both equations
\[
   \sum_{p \pdiv n} \frac{1}{p} \pm \frac{1}{n} = 1 \,.
\]

\section{Bernoulli and Stirling numbers}
\setcounter{equation}{0}

First, we will introduce the Stirling numbers of the first and second kind,
whereas we only need the latter numbers which are connected with Bernoulli numbers.
At the end of this section we will prove Theorem \ref{theor-snm-mbn}.

\begin{defin}
Define falling factorials by
\[
    (x)_n := x(x-1)\cdots(x-n+1) \,, \qquad (x)_0 := 1 \,, \qquad n \in \NN \,.
\]
The Stirling numbers $\stirling_1$ of the first kind and
$\stirling_2$ of the second kind are defined by
\begin{equation} \label{def-stirl}
  (x)_n = \sum_{k=0}^n \stirling_1 (n,k) \, x^k \,, \qquad
  x^n = \sum_{k=0}^n \stirling_2 (n,k) \, (x)_k \,.
\end{equation}
\end{defin}
We basically have with $\nu=1, 2$
\begin{equation} \label{eqn-stirl-basic}
   \stirling_\nu (n,k) = \left\{
     \begin{array}{rl}
       1 \,, & n=k \geq 0 \,, \\
       0 \,, & k>n \quad \mbox{or} \quad k=0, n \geq 1 \,.
     \end{array} \right.
\end{equation}

We use the notations $\{\ \}$ und $[\ \ ]$ like in \cite{graham94concrete},
extending to $\langle\ \ \rangle$
\begin{equation} \label{eqn-stirl-nota-def}
   \ste{n}{k} := \stirling_1 (n,k) \,, \qquad
   \stz{n}{k} := \stirling_2 (n,k) \qquad \mbox{and} \qquad \ST{n}{k} := k!
     \, \stz{n}{k} \,.
\end{equation}
By definition of the binomial coefficients, regarded as polynomials
\[
    \binom{x}{k} = \frac{x(x-1)\cdots(x-k+1)}{k!} = \frac{(x)_k}{k!} \,,
\]
equations \refeqn{def-stirl} now become
\begin{equation} \label{eqn-stirl-binom}
  \binom{x}{n} = \sum_{k=0}^n \frac{1}{n!} \ste{n}{k} \, x^k \,, \qquad
    x^n = \sum_{k=0}^n \ST{n}{k} \binom{x}{k} \,.
\end{equation}

Recall the definition of the summation function $S_n$. Now, this function
can be directly derived by Stirling numbers of the second kind.
Similarly, one can deduce an iterated summation function $S_{n,r}$ summing over
$S_{n,{r-1}}$ with $S_{n,1} = S_n$ and $S_{n,0}(m) = m^n$.

\begin{theorem} \label{theor-form-sn-bn}
Let $n$ be a positive integer, then for real $x$
\[
   S_n(x) = \sum_{k=1}^n \ST{n}{k} \binom{x}{k+1}  \,, \quad
   B_n    = \sum_{k=1}^n \ST{n}{k} \frac{(-1)^k}{k+1} \,.
\]
\end{theorem}

\begin{proof}
Let $m > 1$ be a positive integer, then by \refeqn{eqn-stirl-binom} and
summarizing binomial coefficients via Pascal's triangle, we obtain
\[
   S_n(m) = \sum_{\nu=1}^{m-1} \nu^n = \sum_{\nu=1}^{m-1} \sum_{k=1}^n \ST{n}{k} \binom{\nu}{k}
          = \sum_{k=1}^n \ST{n}{k} \sum_{\nu=1}^{m-1} \binom{\nu}{k}
          = \sum_{k=1}^n \ST{n}{k} \binom{m}{k+1} \,.
\]
This shows that $S_n$ is a polynomial of degree $n+1$ which
is then also valid for real $x$ and equals the former function $S_n$ in
Theorem \ref{theor-sn-bn}. Since $S_n(0) = 0$ and $\binom{-1}{k}=(-1)^k$, we can
write
\[
   S_n'(0) = \lim_{x \to 0} \frac{S_n(x)}{x}
           = \lim_{x \to 0} \sum_{k=1}^n \ST{n}{k} \frac{1}{k+1} \binom{x-1}{k}
           = \sum_{k=1}^n \ST{n}{k} \frac{(-1)^k}{k+1} \,.
\]
On the other side, by Theorem \ref{theor-sn-bn} and basic properties of Bernoulli polynomials
\[
   S_n'(0) = \frac{d}{dx} \left( \frac{B_{n+1}(x) - B_{n+1}}{n+1} \right) \!\!\! \valueat{x=0}
           = B_n(0) = B_n \,. \vspace*{-4ex}
\]
\end{proof}

\begin{remark}
The shortest proof of the formula of $B_n$ is given by $p$-adic theory.
For definitions see Robert \cite[Chapter 4/5]{robert00padic}.
By Volkenborn integral, Mahler series, and Stirling numbers,
it immediately follows for any prime $p$ and $n \geq 1$
\[
   B_n = \int_{\ZZ_p} x^n \, dx
       = \int_{\ZZ_p} \left( \sum_{k=1}^n \ST{n}{k} \binom{x}{k} \right) dx
       = \sum_{k=1}^n \ST{n}{k} \frac{(-1)^k}{k+1} \,.
\]
Knowing that for $n, k \geq 1$, see \cite{graham94concrete},
\begin{equation} \label{eqn-stirl-sum}
   \ST{n}{k} = \sum_{\nu = 1}^k \binom{k}{\nu} (-1)^{k-\nu} \, \nu^n \,,
\end{equation}
one can derive a double sum
\[
   B_n = \sum_{k=1}^n \frac{1}{k+1} \,
           \sum_{\nu = 1}^k \binom{k}{\nu} (-1)^{\nu} \, \nu^n
\]
which was already given by Worpitzky \cite{worp83bernoulli} in 1883.
\end{remark}

\begin{lemma} \label{lem-binom-rp-1}
Let $r$ be a positive integer and $p$ be a prime. Let $0 \leq \nu < p$, then
\[
   \binom{rp-1}{\nu} \equiv (-1)^\nu \pmod{p} \,.
\]
\end{lemma}

\begin{proof}
Case $\nu=0$ is trivial. Since $p \notdiv \nu!$
\[
   \binom{rp-1}{\nu} = \frac{(rp-1)\cdots(rp-\nu)}{\nu!}
      \equiv (-1)^{\nu} \, \frac{1 \cdots \nu}{\nu!}
      \equiv (-1)^{\nu} \pmod{p} \,.  \vspace*{-4ex}
\]
\end{proof}

\begin{lemma} \label{lem-stirl-congr}
Let $n, k$ be positive integers with even $n$. Then
\[
     \ST{n}{k-1} \equiv \left\{
     \begin{array}{rl}
       -1 \,, & k=p, \, p-1 \pdiv n \\
       0  \,, & \mbox{otherwise}
     \end{array} \right.
       \pmod{k} \,.
\]
\end{lemma}

\begin{proof}
Consider \refeqn{eqn-stirl-basic}, \refeqn{eqn-stirl-nota-def}, and \refeqn{eqn-stirl-sum}.
Cases $k=1, 2$ are trivial. It is well-known that
\[
   (k-1)! \equiv \left\{
     \begin{array}{rl}
       -1 \,, & k=p \\
       2  \,, & k=4 \\
       0  \,, & \mbox{otherwise}
     \end{array} \right.
       \pmod{k} \,.
\]
Since $(k-1)! \pdiv \STL{n}{k-1}$, cases $k=4$ and
$k=p > 2$ remain. Let $n$ be even with $n \geq 2$. \par
Case $k=4$: \refeqn{eqn-stirl-sum} yields
\[
   \ST{n}{3} = \binom{3}{1} 1^n - \binom{3}{2} 2^n + \binom{3}{3} 3^n
       = 3 - 3 \cdot 2^n + 3^n \equiv -1 + (-1)^n \equiv 0 \pmod{4} \,.
\]
Case $k=p > 2$: Lemma \ref{lem-binom-rp-1} and \refeqn{eqn-stirl-sum} provide
\[
   \ST{n}{p-1} \equiv \sum_{\nu = 1}^{p-1} \binom{p-1}{\nu} (-1)^{\nu} \, \nu^n
     \equiv \sum_{\nu = 1}^{p-1} \nu^n \equiv S_n(p) \equiv
     \left\{ \!\!
        \begin{array}{rl}
          -1 \,, & p-1 \pdiv n \\
          0  \,, & p-1 \notdiv n
        \end{array} \right.
          \!\!\! \pmod{p} \,.
\]
The last part of the congruence is easily derived,
see \cite[Lemma 2, p.\,235]{ireland90classical}.
\end{proof}

Now, we are ready to prove Theorem \ref{theor-snm-mbn}. Note that we do not need
the Theorem of Clausen-von Staudt which will then follow as a corollary.

\begin{proofof}{Theorem \ref{theor-snm-mbn}}
Let $n$, $m$ be positive integers with $m > 1$ and even $n$. We have to show that
\begin{equation} \label{eqn-snm-mbn-sum}
   S_n(m) \equiv m \, B_n \equiv
     - \!\! \sum_{\ds p \pdiv m \atop \ds p-1 \pdiv n} \!\! \frac{m}{p} \pmod{m} \,.
\end{equation}
By Theorem \ref{theor-form-sn-bn} we can write
\[
   S_n(m) = \sum_{k=1}^n \ST{n}{k} \binom{m}{k+1}
     = \sum_{k=2}^{n+1} \ST{n}{k-1} \frac{m}{k} \binom{m-1}{k-1} \,.
\]
If $(k,m)=1$ then $m/k \equiv 0 \ \mod{m}$. Lemma \ref{lem-stirl-congr} states
$k \pdiv \STL{n}{k-1}$ for all $k$ but $(*)$ $k=p$ with $p-1 \pdiv n$.
For all these $k$ except $(*)$ it follows
\begin{equation} \label{eqn-loc-snm-mbn-1}
   \ST{n}{k-1} \frac{m}{k} \binom{m-1}{k-1} \equiv 0 \pmod{m} \,.
\end{equation}
Hence, the following terms remain
\[
   S_n(m) \equiv \sum_{\ds p \pdiv m \atop \ds p-1 \pdiv n}
     \ST{n}{p-1} \frac{m}{p} \binom{m-1}{p-1} \pmod{m} \,.
\]
Now, we use Lemma \ref{lem-divcongr} to evaluate the congruence above.
Since $p \pdiv m$ and $p-1 \pdiv n$, Lemma \ref{lem-binom-rp-1} and
Lemma \ref{lem-stirl-congr} provide, also valid in case $p=2$,
\begin{equation} \label{eqn-loc-snm-mbn-2}
   \ST{n}{p-1} \binom{m-1}{p-1} \equiv (-1) \cdot (-1)^{p-1} \equiv -1 \pmod{p} \,.
\end{equation}
Finally, we obtain
\[
   S_n(m) \equiv \sum_{\ds p \pdiv m \atop \ds p-1 \pdiv n}
     \ST{n}{p-1} \frac{m}{p} \binom{m-1}{p-1}
     \equiv - \!\! \sum_{\ds p \pdiv m \atop \ds p-1 \pdiv n} \frac{m}{p}
       \pmod{m} \,.
\]
On the other side, we can use similar arguments regarding \refeqn{eqn-loc-snm-mbn-1}
and \refeqn{eqn-loc-snm-mbn-2}
\[
   m B_n \equiv \sum_{k=2}^{n+1} \ST{n}{k-1} \frac{m}{k} (-1)^{k-1}
     \equiv \!\! \sum_{\ds p \pdiv m \atop \ds p-1 \pdiv n} \ST{n}{p-1} \frac{m}{p} (-1)^{p-1}
     \equiv - \!\! \sum_{\ds p \pdiv m \atop \ds p-1 \pdiv n} \frac{m}{p} \pmod{m} \,.
      \vspace*{-4ex}
\]
\end{proofof}

\begin{corl}[Clausen-von Staudt]
Let $n$ be an even positive integer. Then
\[
   B_n + \sum_{p-1 \pdiv n} \frac{1}{p} \in \ZZ \qquad \mbox{and} \qquad
   \denom(B_n) = \prod_{p-1 \pdiv n} p \,.
\]
\end{corl}

\begin{proof}
Since \refeqn{eqn-snm-mbn-sum} is valid for arbitrary $m > 1$,
it is also valid for any prime $p$ with $m=p$. The congruence then is $p$-integral
\[
   p B_n \equiv \left\{ \!\!
        \begin{array}{rl}
          -1 \,, & p-1 \pdiv n \\
          0  \,, & p-1 \notdiv n
        \end{array} \right.
          \!\!\! \pmod{p}
\]
which shows that the denominator of $B_n$ must be squarefree
and has the form claimed above.
Now, take $m=\denom(B_n)$ then \refeqn{eqn-snm-mbn-sum} yields, respectively divided by $m$
\[
   m B_n \equiv - \!\! \sum_{p-1 \pdiv n} \frac{m}{p} \pmod{m} \,,
     \quad B_n \equiv - \!\! \sum_{p-1 \pdiv n} \frac{1}{p} \pmod{\ZZ} \,. \vspace*{-4ex}
\]
\end{proof}

\begin{remark} Let $n$ be an even integer and $B_n=U_n/V_n$ with $(U_n,V_n)=1$ and $V_n > 0$, then
the last congruence reads
\[
   U_n \equiv - \!\! \sum_{p-1 \pdiv n} \frac{V_n}{p} \pmod{V_n} \,.
\]
\end{remark}

\subsection*{Remark}

This article is based on a part of
the author's previous diploma thesis \cite[Chapter 2/3]{kellner02irrpairord},
where some results are described more generally.
\bigskip

Bernd C. Kellner \\
address: Reitstallstr. 7, 37073 G\"ottingen, Germany \\
email: bk@bernoulli.org

\bibliographystyle{alpha}
\bibliography{eq_bib}

\end{document}